# Dimension and product structure of hyperbolic measures

By Luis Barreira, Yakov Pesin, and Jörg Schmeling*

**Abstract**

We prove that every hyperbolic measure invariant under a $C^{1+\alpha}$ diffeomorphism of a smooth Riemannian manifold possesses asymptotically "almost" local product structure, i.e., its density can be approximated by the product of the densities on stable and unstable manifolds up to small exponentials. This has not been known even for measures supported on locally maximal hyperbolic sets.

Using this property of hyperbolic measures we prove the long-standing Eckmann-Ruelle conjecture in dimension theory of smooth dynamical systems: the pointwise dimension of every hyperbolic measure invariant under a $C^{1+\alpha}$ diffeomorphism exists almost everywhere. This implies the crucial fact that virtually all the characteristics of dimension type of the measure (including the Hausdorff dimension, box dimension, and information dimension) coincide. This provides the rigorous mathematical justification of the concept of fractal dimension for hyperbolic measures.

## 1. Introduction

In this paper we provide an affirmative solution of the long-standing problem in the interface of dimension theory and dynamical systems known as the Eckmann-Ruelle conjecture.

In the late 70's–beginning 80's, attention of many physicists and applied mathematicians had turned to the study of dimension of strange attractors

*This paper was written while L. B. was on leave from Instituto Superior Técnico, Department of Mathematics, at Lisbon, Portugal, and J. S. was visiting Penn State. L. B. was partially supported by FCT's Pluriannual Funding Program and PRAXIS XXI grants 2/2.1/MAT/199/94, BD5236/95, and PBIC/C/MAT/2139/95. Ya. P. was partially supported by the National Science Foundation grant #DMS9403723. J. S. was supported by the Leopoldina-Forderpreis. L. B. and Ya. P. were partially supported by the NATO grant CRG970161.

1991 *Mathematics Subject Classification*. 58F11, 28D05.

(i.e., attracting invariant sets with some hyperbolic structure) in evolution-type systems (see, for example, [6], [7], [22]). The dimension was used to characterize a (finite) number of independent modes needed to describe the infinite-dimensional system. Several results were obtained which indicated relations between the dimension of the attractor and other invariants of the system (such as Lyapunov exponents and entropy; see, for example, [13], [14], [16], [23]). This study has become an important breakthrough in understanding the structure of systems of evolution type.

In the survey article [4], Eckmann and Ruelle summarized this activity and outlined a rigorous mathematical foundation for it. They considered dynamical systems with chaotic behavior of trajectories and described relations between persistence of chaotic motions and existence of strange attractors. They also discussed various concepts of dimension and pointed out the importance of the so-called *pointwise (local) dimension* of invariant measures. For a Borel measure $\mu$ on a complete metric space $M$, the latter is defined by

$$(1) \qquad d(x) \stackrel{\text{def}}{=} \lim_{r \to 0} \frac{\log \mu(B(x,r))}{\log r}$$

where $B(x,r)$ is the ball centered at $x$ of radius $r$ (provided the limit exists). It was introduced by Young in [23] and characterizes the local geometrical structure of an invariant measure with respect to the metric in the phase space of the system. Its crucial role in dimension theory of dynamical systems was acknowledged by many experts in the field (see, for example, the paper by Farmer, Ott, and Yorke [6], the ICM address by Young [24, p. 1232] and also [25, p. 318]).

If the limit in (1) does not exist one can consider the lower and upper limits and introduce respectively *the lower and upper pointwise dimensions* of $\mu$ at $x$ which we denote by $\underline{d}(x)$ and $\overline{d}(x)$.

The existence of the limit in (1) for a Borel probability measure $\mu$ on $M$ implies the crucial fact that virtually *all* the known characteristics of dimension type of the measure coincide (this is partly described in Prop. 1 in §2). The common value is a fundamental characteristic of the fractal structure of $\mu$ — the *fractal dimension* of $\mu$.

In this paper we consider a $C^{1+\alpha}$ diffeomorphism of a compact smooth Riemannian manifold without boundary. Our goal is to show the existence of the pointwise dimension in the case when $\mu$ is *hyperbolic*, i.e., all the Lyapunov exponents of $f$ are nonzero at $\mu$-almost every point (see Main Theorem in §3). This statement has been an open problem in dimension theory of dynamical systems for about 15 years and is often referred to as the *Eckmann-Ruelle conjecture*.

Since hyperbolic measures play a crucial role in studying physical models with persistent chaotic behavior and fractal structure of invariant sets, our result provides a rigorous mathematical foundation for such a study.

The problem of the existence of the pointwise dimension has a long history. In [23], Young obtained a positive answer for a hyperbolic measure $\mu$ invariant under a $C^{1+\alpha}$ surface diffeomorphism $f$. Moreover, she showed that in this case for almost every point $x$,

$$\underline{d}(x) = \overline{d}(x) = h_\mu(f) \left( \frac{1}{\lambda_1} - \frac{1}{\lambda_2} \right),$$

where $h_\mu(f)$ is the metric entropy of $f$ and $\lambda_1 > 0 > \lambda_2$ are the Lyapunov exponents of $\mu$.

In [12], Ledrappier established the existence of the pointwise dimension for arbitrary SRB-*measures* (called so after Sinai, Ruelle, and Bowen). In [20], Pesin and Yue extended his approach and proved the existence for hyperbolic measures satisfying the so-called *semi-local product structure* (this class includes, for example, *Gibbs measures* on locally maximal hyperbolic sets).

A substantial breakthrough in studying the pointwise dimension was made by Ledrappier and Young in [16]. They proved the existence of the *stable* and *unstable pointwise dimensions*, i.e., the pointwise dimensions along stable and unstable local manifolds for typical points (see Prop. 2 in §2). They also showed that the upper pointwise dimension at a typical point does not exceed the sum of the stable and unstable pointwise dimensions.

Our proof exploits their result in an essential way. It also uses a new and nontrivial property of hyperbolic ergodic measures that we establish in this paper. Loosely speaking, this property means that such measures have asymptotically "almost" local product structure. Let us point out that this property has not been known even for invariant measures on locally maximal hyperbolic sets (whose local topological structure is the direct product). This property also enables us to show that the pointwise dimension of a hyperbolic measure is almost everywhere the sum of the pointwise dimensions along stable and unstable local manifolds.

*Acknowledgment.* We would like to thank François Ledrappier for useful discussions and comments.

## 2. Preliminaries

2.1. *Facts from dimension theory.* We describe some most important characteristics of dimension type (see, for example, [5], [18]). Let $X$ be a complete separable metric space. For a subset $Z \subset X$ and a number $\alpha \geq 0$ the $\alpha$-*Hausdorff measure of* $Z$ is defined by

$$m_H(Z, \alpha) = \liminf_{\varepsilon \to 0} \sum_{\mathcal{G}} (\operatorname{diam} U)^\alpha,$$

where the infimum is taken over all finite or countable coverings $\mathcal{G}$ of $Z$ by open sets with $\operatorname{diam} \mathcal{G} \leq \varepsilon$. The *Hausdorff dimension of $Z$* (denoted $\dim_H Z$) is defined by

$$\dim_H Z = \inf\{\alpha : m_H(Z, \alpha) = 0\} = \sup\{\alpha : m_H(Z, \alpha) = \infty\}.$$

We define the *lower* and *upper box dimensions* of $Z$ (denoted respectively by $\underline{\dim}_B Z$ and $\overline{\dim}_B Z$) by

$$\underline{\dim}_B Z = \inf\{\alpha : \underline{r}_H(Z, \alpha) = 0\} = \sup\{\alpha : \underline{r}_H(Z, \alpha) = \infty\},$$

$$\overline{\dim}_B Z = \inf\{\alpha : \overline{r}_H(Z, \alpha) = 0\} = \sup\{\alpha : \overline{r}_H(Z, \alpha) = \infty\},$$

where

$$\underline{r}_H(Z, \alpha) = \liminf_{\varepsilon \to 0} \sum_{U \in \mathcal{G}} \varepsilon^\alpha, \qquad \overline{r}_H(Z, \alpha) = \overline{\liminf_{\varepsilon \to 0}} \sum_{U \in \mathcal{G}} \varepsilon^\alpha$$

and the infimum is taken over all finite or countable coverings $\mathcal{G}$ of $Z$ by open sets of diameter $\varepsilon$. One can show that

$$\underline{\dim}_B Z = \lim_{\varepsilon \to 0} \frac{\log N(Z, \varepsilon)}{\log(1/\varepsilon)}, \quad \overline{\dim}_B Z = \overline{\lim_{\varepsilon \to 0}} \frac{\log N(Z, \varepsilon)}{\log(1/\varepsilon)},$$

where $N(Z, \varepsilon)$ is the smallest number of balls of radius $\varepsilon$ needed to cover the set $Z$.

It is easy to see that

$$\dim_H Z \leq \underline{\dim}_B Z \leq \overline{\dim}_B Z.$$

The coincidence of the Hausdorff dimension and lower and upper box dimension is a relatively rare phenomenon and can occur only in some "rigid" situations (see [1], [5], [19]).

In order to describe the geometric structure of a subset $Z$ invariant under a dynamical system $f$ acting on $X$ we consider a measure $\mu$ supported on $Z$. Its *Hausdorff dimension* and *lower* and *upper box dimensions* (which are denoted by $\dim_H \mu$, $\underline{\dim}_B \mu$, and $\overline{\dim}_B \mu$, respectively) are

$$\dim_H \mu = \inf\{\dim_H Z : \mu(Z) = 1\},$$

$$\underline{\dim}_B \mu = \lim_{\delta \to 0} \inf\{\underline{\dim}_B Z : \mu(Z) \geq 1 - \delta\},$$

$$\overline{\dim}_B \mu = \lim_{\delta \to 0} \inf\{\overline{\dim}_B Z : \mu(Z) \geq 1 - \delta\}.$$

From the definition it follows that

$$\dim_H \mu \leq \underline{\dim}_B \mu \leq \overline{\dim}_B \mu.$$

Another important characteristic of dimension type of $\mu$ is its *information dimension*. Given a partition $\xi$ of $X$, define the *entropy of $\xi$ with respect to $\mu$* by

$$H_\mu(\xi) = -\sum_{C_\xi} \mu(C_\xi) \log \mu(C_\xi)$$

where $C_\xi$ is an element of the partition $\xi$. Given a number $\varepsilon > 0$, set

$$H_\mu(\varepsilon) = \inf\{H_\mu(\xi) : \operatorname{diam} \xi \leq \varepsilon\}$$

where $\operatorname{diam} \xi = \max \operatorname{diam} C_\xi$. We define the *lower* and *upper information dimensions of $\mu$* by

$$\underline{I}(\mu) = \varliminf_{\varepsilon \to 0} \frac{H_\mu(\varepsilon)}{\log(1/\varepsilon)}, \quad \overline{I}(\mu) = \varlimsup_{\varepsilon \to 0} \frac{H_\mu(\varepsilon)}{\log(1/\varepsilon)}.$$

There is a powerful criterion established by Young in [23] that guarantees the coincidence of the Hausdorff dimension and lower and upper box dimensions of measures as well as their lower and upper information dimensions.

PROPOSITION 1 ([23]). *Let $X$ be a compact separable metric space of finite topological dimension and $\mu$ a Borel probability measure on $X$. Assume that*

(2) $$\underline{d}(x) = \overline{d}(x) = d$$

*for $\mu$-almost every $x \in X$. Then*

$$\dim_H \mu = \underline{\dim}_B \mu = \overline{\dim}_B \mu = \underline{I}(\mu) = \overline{I}(\mu) = d.$$

A measure $\mu$ which satisfies (2) is called *exact dimensional*.

2.2. *Hyperbolic measures.* Let $M$ be a smooth Riemannian manifold without boundary, and $f \colon M \to M$ a $C^{1+\alpha}$ diffeomorphism on $M$. Let also $\mu$ be an $f$-invariant ergodic Borel probability measure on $M$.

Given $x \in M$ and $v \in T_x M$ define the *Lyapunov exponent of $v$ at $x$* by the formula

$$\lambda(x, v) = \varlimsup_{n \to \infty} \frac{1}{n} \log \|d_x f^n v\|.$$

If $x$ is fixed then the function $\lambda(x, \cdot)$ can take on only finitely many values $\lambda_1(x) < \cdots < \lambda_{k(x)}(x)$. The functions $\lambda_i(x)$ are measurable and $f$-invariant. Since $\mu$ is ergodic, these functions are constant $\mu$-almost everywhere. We denote these constants by $\lambda_1 < \cdots < \lambda_k$. The measure $\mu$ is said to be *hyperbolic* if $\lambda_i \neq 0$ for every $i = 1, \ldots, k$.

There exists a measurable function $r(x) > 0$ such that for $\mu$-almost every $x \in M$ the sets

$$W^s(x) = \left\{ y \in B(x, r(x)) : \varlimsup_{n \to +\infty} \frac{1}{n} \log d(f^n x, f^n y) < 0 \right\},$$

$$W^u(x) = \left\{ y \in B(x, r(x)) : \varlimsup_{n \to -\infty} \frac{1}{n} \log d(f^n x, f^n y) > 0 \right\}$$

are immersed local manifolds called *stable* and *unstable local manifolds* at $x$ (see [17] for details). For each $r \in (0, r(x))$ we consider the balls $B^s(x, r) \subset W^s(x)$ and $B^u(x, r) \subset W^u(x)$ centered at $x$ with respect to the induced distances on $W^s(x)$ and $W^u(x)$ respectively.

Let $\xi$ be a measurable partition of $M$. It has a canonical system of conditional measures: for $\mu$-almost every $x$ there is a probability measure $\mu_x$ defined on the element $\xi(x)$ of $\xi$ containing $x$. The conditional measures $\mu_x$ are uniquely characterized by the following property: if $\mathcal{B}_\xi$ is the $\sigma$-subalgebra (of the Borel $\sigma$-algebra) whose elements are unions of elements of $\xi$, and $A \subset M$ is a measurable set, then $x \mapsto \mu_x(A \cap \xi(x))$ is $\mathcal{B}_\xi$-measurable and

$$\mu(A) = \int \mu_x(A \cap \xi(x)) \, d\mu(x).$$

In [16], Ledrappier and Young constructed two measurable partitions $\xi^s$ and $\xi^u$ of $M$ such that for $\mu$-almost every $x \in M$:

1. $\xi^s(x) \subset W^s(x)$ and $\xi^u(x) \subset W^u(x)$;
2. $\xi^s(x)$ and $\xi^u(x)$ contain the intersection of an open neighborhood of $x$ with $W^s(x)$ and $W^u(x)$ respectively.

We denote the system of conditional measures of $\mu$ with respect to the partitions $\xi^s$ and $\xi^u$, respectively by $\mu_x^s$ and $\mu_x^u$, and for any measurable set $A \subset M$ we write $\mu_x^s(A) = \mu_x^s(A \cap \xi^s(x))$ and $\mu_x^u(A) = \mu_x^u(A \cap \xi^u(x))$.

Given $x \in M$, consider the lower and upper pointwise dimensions of $\mu$ at $x$, $\underline{d}(x)$ and $\overline{d}(x)$. Since these functions are measurable and $f$-invariant they are constant $\mu$-almost everywhere. We denote these constants by $\underline{d}$ and $\overline{d}$ respectively.

In [16], Ledrappier and Young introduced the quantities

$$d^s(x) \stackrel{\text{def}}{=} \lim_{r \to 0} \frac{\log \mu_x^s(B^s(x, r))}{\log r},$$

$$d^u(x) \stackrel{\text{def}}{=} \lim_{r \to 0} \frac{\log \mu_x^u(B^u(x, r))}{\log r}$$

provided that the corresponding limits exist at $x \in M$. We call them, respectively, *stable* and *unstable pointwise dimensions* of $\mu$.

PROPOSITION 2 ([16]).
1. *For $\mu$-almost every $x \in M$ the limits $d^s(x)$ and $d^u(x)$ exist and are constant $\mu$-almost everywhere; we denote these constants by $d^s$ and $d^u$.*
2. *If $\mu$ is a hyperbolic measure then*

$$\overline{d} \leq d^s + d^u.$$

When the entropy of $f$ is zero it follows from [16] that $d^s = d^u = 0$ and hence $\underline{d} = \overline{d} = d^s + d^u = 0$.

Let us point out that in [16] the authors consider a class of measures more general than hyperbolic measures (some of the Lyapunov exponents may be zero). They prove Proposition 2 under the assumption that the diffeomorphism $f$ is of class $C^2$. The main ingredient of their proof is the existence of *intermediate* pointwise dimensions, i.e., the pointwise dimensions of the conditional measures generated by the invariant measure on the intermediate stable and unstable leaves. This, in turn, relies on the Lipschitz continuity of the holonomy map generated by these intermediate leaves (see the Appendix for definitions and precise statements) and is where the assumption that $f$ is of class $C^2$ is used (see §(4.2) in [16]; other arguments in [16] do not use the Lipschitz continuity of the holonomy map and go well for $C^{1+\alpha}$ diffeomorphisms).

In the Appendix to the paper we provide a proof of the fact that for hyperbolic measures the Lipschitz property of the holonomy map holds for diffeomorphisms of class $C^{1+\alpha}$. Indeed, we prove a slightly more general statement: the Lipschitz property holds for intermediate stable and unstable foliations even if some of the Lyapunov exponents are zero. Our approach gives a new proof of this property even in the case of diffeomorphisms of class $C^2$. As a consequence we obtain that Proposition 2 holds for diffeomorphisms of class $C^{1+\alpha}$ and so does our Main Theorem.

Let us point out that the Lipschitz continuity of the holonomy map presumably fails if the map is generated by stable leaves inside the stable-neutral foliation.

## 3. Main Theorem

In this paper we prove the following statement.

MAIN THEOREM. *Let $f$ be a $C^{1+\alpha}$ diffeomorphism on a smooth Riemannian manifold $M$ without boundary, and $\mu$ an $f$-invariant compactly supported ergodic Borel probability measure. If $\mu$ is hyperbolic then the following properties hold:*

1. *for every $\delta > 0$, there exist a set $\Lambda \subset M$ with $\mu(\Lambda) > 1 - \delta$ and a constant $\kappa \geq 1$ such that for every $x \in \Lambda$ and every sufficiently small $r$ (depending on $x$),*

$$
\text{(3)} \quad r^\delta \mu_x^s(B^s(x, \frac{r}{\kappa})) \mu_x^u(B^u(x, \frac{r}{\kappa})) \leq \mu(B(x, r)) \leq r^{-\delta} \mu_x^s(B^s(x, \kappa r)) \mu_x^u(B^u(x, \kappa r));
$$

2. *$\mu$ is exact dimensional and its pointwise dimension is equal to the sum of the stable and unstable pointwise dimensions, i.e.,*

$$\underline{d} = \overline{d} = d^s + d^u.$$

This statement provides an affirmative solution to the Eckmann-Ruelle conjecture and describes the most broad class of measures invariant under smooth dynamical systems which are exact dimensional.

Note that an SRB-measure is locally equivalent on $\Lambda$ to the direct product of an absolutely continuous measure on an unstable leaf and a measure on a stable leaf (see [12]). Hence, statement 1 holds in this case automatically. One can also show that Gibbs measures on a locally maximal hyperbolic set $\Lambda$ are locally equivalent to the direct product of a measure on an unstable leaf and a measure on a stable leaf (see [9]). Therefore, statement 1 holds in this case as well.

Let us also point out that neither of the assumptions of the Main Theorem can be omitted. Ledrappier and Misiurewicz [15] constructed an example of a smooth map of a circle preserving an ergodic measure with zero Lyapunov exponent which is not exact dimensional. In [19], Pesin and Weiss presented an example of a Hölder homeomorphism with Hölder constant arbitrarily close to 1 whose measure of maximal entropy is not exact dimensional.

*Remarks.* 1. Statement 1 of the Main Theorem establishes a new and nontrivial property of an arbitrary hyperbolic measure. Loosely speaking, it means that every hyperbolic invariant measure possesses asymptotically "almost" local product structure. This statement has not been known even for measures supported on (uniformly) hyperbolic locally maximal invariant sets. The lower bound in (3) can be easily obtained from results in [16] while the upper bound is one of the main ingredients of our proof. Note that statement 2 follows from statement 1. The proof of statement 2 exploits the existence of stable and unstable pointwise dimensions and the argument in [20] (see §6).

In order to illustrate the property of having asymptotically "almost" local product structure, let us consider an ergodic measure $\mu$ invariant under the full shift $\sigma$ on the space $\Sigma_p$ of all two-sided infinite sequences of $p$ numbers. This space is endowed with the usual "symbolic" metric $d_\beta$, for each fixed number

$\beta \geq 2$, defined by
$$d_\beta(\omega^1, \omega^2) = \sum_{i \in \mathbb{Z}} \beta^{-|i|} |\omega_i^1 - \omega_i^2|,$$
where $\omega^1 = (\omega_i^1)$ and $\omega^2 = (\omega_i^2)$. Fix $\omega = (\omega_i) \in \Sigma_p$. The cylinder
$$C_n(\omega) = \{\bar\omega = (\bar\omega_i) : \bar\omega_i = \omega_i \text{ for } i = -n, \ldots, n\}$$
can be identified with the direct product $C_n^+(\omega) \times C_n^-(\omega)$ where
$$C_n^+(\omega) = \{\bar\omega = (\bar\omega_i) : \bar\omega_i = \omega_i \text{ for } i = 0, \ldots, n\}$$
and
$$C_n^-(\omega) = \{\bar\omega = (\bar\omega_i) : \bar\omega_i = \omega_i \text{ for } i = -n, \ldots, 0\}$$
are the "positive" and "negative" cylinders at $\omega$ of "size" $n$. Define measures
$$\mu_n^+(\omega) = \mu|C_n^+(\omega) \quad \text{and} \quad \mu_n^-(\omega) = \mu|C_n^-(\omega).$$
The measure $\mu$ is said to have local product structure if the measure $\mu|C_n(\omega)$ is equivalent to the product $\mu_n^+(\omega) \times \mu_n^-(\omega)$ uniformly over $\omega \in \Sigma_p$ and $n > 0$. It is known that Gibbs measures have local product structure (see, for example, [18]). For an arbitrary $\sigma$-invariant ergodic measure $\mu$ it follows from statement 1 of the Main Theorem (see (3)) that for every $\delta > 0$ there exist a set $\Lambda \subset \Sigma_p$ with $\mu(\Lambda) > 1 - \delta$ and an integer $m \geq 1$ such that for every $\omega \in \Lambda$ and every sufficiently large $n$ (depending on $\omega$),
$$\beta^{-\delta|n|} \mu_{n+m}^+(\omega) \times \mu_{n+m}^-(\omega) \leq \mu|C_n(\omega) \leq \beta^{\delta|n|} \mu_{n-m}^+(\omega) \times \mu_{n-m}^-(\omega).$$

2. It follows immediately from the Main Theorem that the pointwise dimension of an ergodic invariant measure supported on a (uniformly) hyperbolic locally maximal invariant set is exact dimensional. This result has not been known before. We emphasize that in this situation the stable and unstable foliations need not be Lipschitz (in fact, they are "generically" not Lipschitz; see [21]), and, in general, the measure need not have a local product structure despite the fact that the set itself does. Therefore, both statements of the Main Theorem are nontrivial even for measures supported on hyperbolic locally maximal invariant sets.

3. The role of the Eckmann-Ruelle conjecture in dimension theory of dynamical systems is similar to the role of the Shannon-McMillan-Breiman theorem in the entropy theory. In order to illustrate this, consider the full shift $\sigma$ on the space $\Sigma_p$.

Let $\mu$ be a $\sigma$-invariant ergodic measure on $\Sigma_p$. By the Shannon-McMillan-Breiman theorem, for $\mu$-almost every $\omega \in \Sigma_p$,

(4) $$\lim_{n \to \infty} -\frac{1}{2n+1} \log \mu(C_n(\omega)) = h$$

where $h = h_\mu(\sigma)$. Since the cylinder $C_n(\omega)$ is the ball (in the symbolic metric $d_\beta$) centered at $\omega$ of radius $c\beta^n$ (for some $c > 0$), the quantity in the right-hand side in (4) divided by $\log \beta$ is the pointwise dimension of $\mu$ at $\omega$. The Shannon-McMillan-Breiman theorem (applied to the shift map) thus claims that the pointwise dimension of $\mu$ exists almost everywhere; it is clearly almost everywhere constant; furthermore, the common value is the measure-theoretic entropy of $\mu$ divided by $\log \beta$.

As an important consequence of this theorem one obtains that various definitions of the entropy (due to Kolmogorov and Sinai [11], Katok [8], Brin and Katok [2], etc.) coincide (see [18] for details).

4. Let $\mu$ be an arbitrary (not necessarily ergodic) invariant measure for a $C^{1+\alpha}$ diffeomorphism $f$ of $M$. The measure $\mu$ is said to be *hyperbolic* if on almost every ergodic component the induced measure is a hyperbolic measure on $M$. The Main Theorem remains true for any $f$-invariant compactly supported hyperbolic Borel probability measure, i.e., the pointwise dimension of such a measure exists almost everywhere (but may not be any longer a constant; see §7).

## 4. Description of a special partition

We use the following notation. Let $\eta$ be a partition. For every integers $k, l \geq 1$, we define the partition $\eta_k^l = \bigvee_{n=-k}^{l} f^{-n}\eta$. We observe that $\eta_k^0(x) \cap \eta_0^l(x) = \eta_k^l(x)$.

From now on we assume that $\mu$ is hyperbolic. In [16], Ledrappier and Young constructed a special countable partition $\mathcal{P}$ of $M$ of finite entropy satisfying the following properties. Given $0 < \varepsilon < 1$, there exists a set $\Gamma \subset M$ of measure $\mu(\Gamma) > 1 - \varepsilon/2$, an integer $n_0 \geq 1$, and a number $C > 1$ such that for every $x \in \Gamma$ and any integer $n \geq n_0$, the following statements hold:

a. For all integers $k, l \geq 1$ we have

$$C^{-1} e^{-(l+k)h - (l+k)\varepsilon} \leq \mu(\mathcal{P}_k^l(x)) \leq C e^{-(l+k)h + (l+k)\varepsilon}, \tag{5}$$

$$C^{-1} e^{-kh - k\varepsilon} \leq \mu_x^s(\mathcal{P}_k^0(x)) \leq C e^{-kh + k\varepsilon}, \tag{6}$$

$$C^{-1} e^{-lh - l\varepsilon} \leq \mu_x^u(\mathcal{P}_0^l(x)) \leq C e^{-lh + l\varepsilon}, \tag{7}$$

where $h$ is the Kolmogorov-Sinai entropy of $f$ with respect to $\mu$.

b.

$$\xi^s(x) \cap \bigcap_{n \geq 0} \mathcal{P}_0^n(x) \supset B^s(x, e^{-n_0}), \tag{8}$$

$$\xi^u(x) \cap \bigcap_{n \geq 0} \mathcal{P}_n^0(x) \supset B^u(x, e^{-n_0}). \tag{9}$$

c.

(10) $$e^{-d^s n - n\varepsilon} \leq \mu_x^s(B^s(x, e^{-n})) \leq e^{-d^s n + n\varepsilon},$$

(11) $$e^{-d^u n - n\varepsilon} \leq \mu_x^u(B^u(x, e^{-n})) \leq e^{-d^u n + n\varepsilon}.$$

d.

(12) $$\mathcal{P}_{an}^{an}(x) \subset B(x, e^{-n}) \subset \mathcal{P}(x),$$

(13) $$\mathcal{P}_{an}^{0}(x) \cap \xi^s(x) \subset B^s(x, e^{-n}) \subset \mathcal{P}(x) \cap \xi^s(x),$$

(14) $$\mathcal{P}_{0}^{an}(x) \cap \xi^u(x) \subset B^u(x, e^{-n}) \subset \mathcal{P}(x) \cap \xi^u(x),$$

where $a$ is the integer part of $2(1+\varepsilon)\max\{-\lambda_1, \lambda_k, 1\}$.

e. Define $Q_n(x)$ by

(15) $$Q_n(x) = \bigcup \mathcal{P}_{an}^{an}(y)$$

where the union is taken over all $y \in \Gamma$ for which

$$\mathcal{P}_{0}^{an}(y) \cap B^u(x, 2e^{-n}) \neq \emptyset \text{ and } \mathcal{P}_{an}^{0}(y) \cap B^s(x, 2e^{-n}) \neq \emptyset;$$

then

(16) $$B(x, e^{-n}) \cap \Gamma \subset Q_n(x) \subset B(x, 4e^{-n}),$$

and for each $y \in Q_n(x)$,

$$\mathcal{P}_{an}^{an}(y) \subset Q_n(x).$$

Increasing $n_0$ if necessary we may also assume that

f. For every $x \in \Gamma$ and $n \geq n_0$,

(17) $$B^s(x, e^{-n}) \cap \Gamma \subset Q_n(x) \cap \xi^s(x) \subset B^s(x, 4e^{-n}),$$

(18) $$B^u(x, e^{-n}) \cap \Gamma \subset Q_n(x) \cap \xi^u(x) \subset B^u(x, 4e^{-n}).$$

The above statements are slightly different versions of statements in [16]. Property (5) essentially follows from the Shannon-McMillan-Breiman theorem applied to the partition $\mathcal{P}$ while properties (6) and (7) follow from "leaf-wise" versions of this theorem. The inequalities in (10) and (11) are easy consequences of the existence of the stable and unstable pointwise dimensions $d^s$ and $d^u$ (see Prop. 2). Since the Lyapunov exponents at $\mu$-almost every point are constant and equal to $\lambda_1, \ldots, \lambda_k$, the properties (12), (13), and (14) follow from (8), (9), and the choice of $a$ indicated above. The inclusions in (16) are based upon the continuous dependence of stable and unstable manifolds in the $C^{1+\alpha}$ topology on the base point (in each Pesin set). We need the following well-known result.

PROPOSITION 3 (Borel density lemma). *Let $\mu$ be a finite Borel measure and $A \subset M$ a measurable set. Then for $\mu$-almost every $x \in A$,*

$$\lim_{r \to 0} \frac{\mu(B(x,r) \cap A)}{\mu(B(x,r))} = 1.$$

*Furthermore, if $\mu(A) > 0$ then, for each $\delta > 0$, there is a set $\Delta \subset A$ with $\mu(\Delta) > \mu(A) - \delta$, and a number $r_0 > 0$ such that for all $x \in \Delta$ and $0 < r < r_0$, we have*

$$\mu(B(x,r) \cap A) \geq \frac{1}{2}\mu(B(x,r)).$$

It immediately follows from the Borel density lemma that one can choose an integer $n_1 \geq n_0$ and a set $\widehat{\Gamma} \subset \Gamma$ of measure $\mu(\widehat{\Gamma}) > 1 - \varepsilon$ such that for every $n \geq n_1$ and $x \in \widehat{\Gamma}$,

$$\mu(B(x, e^{-n}) \cap \Gamma) \geq \frac{1}{2}\mu(B(x, e^{-n})); \tag{19}$$

$$\mu_x^s(B^s(x, e^{-n}) \cap \Gamma) \geq \frac{1}{2}\mu_x^s(B^s(x, e^{-n})); \tag{20}$$

$$\mu_x^u(B^u(x, e^{-n}) \cap \Gamma) \geq \frac{1}{2}\mu_x^u(B^u(x, e^{-n})). \tag{21}$$

We establish two additional properties of the partitions $\mathcal{P}_0^k$ and $\mathcal{P}_k^0$.

PROPOSITION 4. *There exists a positive constant $D = D(\widehat{\Gamma}) < 1$ such that for every $k \geq 1$ and $x \in \Gamma$,*

$$\mu_x^s(\mathcal{P}_0^k(x) \cap \Gamma) \geq D;$$

$$\mu_x^u(\mathcal{P}_k^0(x) \cap \Gamma) \geq D.$$

*Proof.* By (8), for every $k \geq 1$ and $x \in \Gamma$, the set $\mathcal{P}_0^k(x) \cap \Gamma$ contains the set $B^s(x, e^{-n_0}) \cap \Gamma$. It follows from (20) and (10) that

$$\mu_x^s(\mathcal{P}_0^k(x) \cap \Gamma) \geq \frac{1}{2}\mu_x^s(B^s(x, e^{-n_0})) \geq \frac{1}{2}e^{-d^s n_0 - n_0 \varepsilon} \stackrel{\text{def}}{=} D.$$

The second inequality in the proposition can be proved in a similar fashion using the properties (9), (21), and (11). $\square$

The next statement establishes the property of the partition $\mathcal{P}$ which simulates the well-known Markov property.

PROPOSITION 5. *For every $x \in \Gamma$ and $n \geq n_0$,*

$$\mathcal{P}_{an}^{an}(x) \cap \xi^s(x) = \mathcal{P}_{an}^0(x) \cap \xi^s(x);$$

$$\mathcal{P}_{an}^{an}(x) \cap \xi^u(x) = \mathcal{P}_0^{an}(x) \cap \xi^u(x).$$

*Proof.* It follows from (13) and (8) that

$$\mathcal{P}^0_{an}(x) \cap \xi^s(x) \subset \mathcal{P}^0_{an}(x) \cap B^s(x, e^{-n}) \subset \mathcal{P}^0_{an}(x) \cap B^s(x, e^{-n_0})$$
$$\subset \mathcal{P}^0_{an}(x) \cap \mathcal{P}^{an}_0(x) \cap \xi^s(x) = \mathcal{P}^{an}_{an}(x) \cap \xi^s(x).$$

Since $\mathcal{P}^{an}_{an}(x) \subset \mathcal{P}^0_{an}(x)$, this completes the proof of the first identity. The proof of the other identity is similar. □

## 5. Preparatory lemmata

Fix $x \in \widehat{\Gamma}$ and an integer $n \geq n_1$. We consider the following two classes $\mathcal{R}(n)$ and $\mathcal{F}(n)$ of elements of the partition $\mathcal{P}^{an}_{an}$ (we call these elements "rectangles"):

$$\mathcal{R}(n) = \{\mathcal{P}^{an}_{an}(y) \subset \mathcal{P}(x) : \mathcal{P}^{an}_{an}(y) \cap \Gamma \neq \emptyset\};$$
$$\mathcal{F}(n) = \{\mathcal{P}^{an}_{an}(y) \subset \mathcal{P}(x) : \mathcal{P}^0_{an}(y) \cap \widehat{\Gamma} \neq \emptyset \text{ and } \mathcal{P}^{an}_0(y) \cap \widehat{\Gamma} \neq \emptyset\}.$$

The rectangles in $\mathcal{R}(n)$ carry all the measure of the set $\mathcal{P}(x) \cap \Gamma$, i.e.,

$$\sum_{R \in \mathcal{R}(n)} \mu(R \cap \Gamma) = \mu(\mathcal{P}(x) \cap \Gamma).$$

Obviously, the rectangles in $\mathcal{R}(n)$ that intersect $\widehat{\Gamma}$ belong to $\mathcal{F}(n)$. If these were the only ones in $\mathcal{F}(n)$, the measure $\mu|\mathcal{P}(x) \cap \Gamma$ would have the "direct product structure" at the "level" $n$. One could then use the approach in [12], [20] to estimate the measure of a ball by the product of its stable and unstable measures. In the general case, the rectangles in the class $\mathcal{F}(n)$ are obtained from the rectangles in $\mathcal{R}(n)$ (that intersect $\widehat{\Gamma}$) by "filling in" the gaps in the "product structure" (see Fig. 1).

We wish to compare the number of rectangles in $\mathcal{R}(n)$ and $\mathcal{F}(n)$ intersecting a given set. This will allow us to evaluate the deviation of the measure $\mu$ from the direct product structure at the level $n$. Our main observation is that for "typical" points $y \in \widehat{\Gamma}$ the number of rectangles from the class $\mathcal{R}(n)$ intersecting $W^s(y)$ (respectively $W^u(y)$) is "asymptotically" the same up to a factor that grows at most subexponentially with $n$.

However, in general, the distribution of these rectangles along $W^s(y)$ (respectively $W^u(y)$) may be different for different points $y$. This causes a deviation from the direct product structure. We will use a simple combinatorial argument to show that this deviation grows at most subexponentially with $n$. One can then say that the measure $\mu$ has an "almost direct product structure."

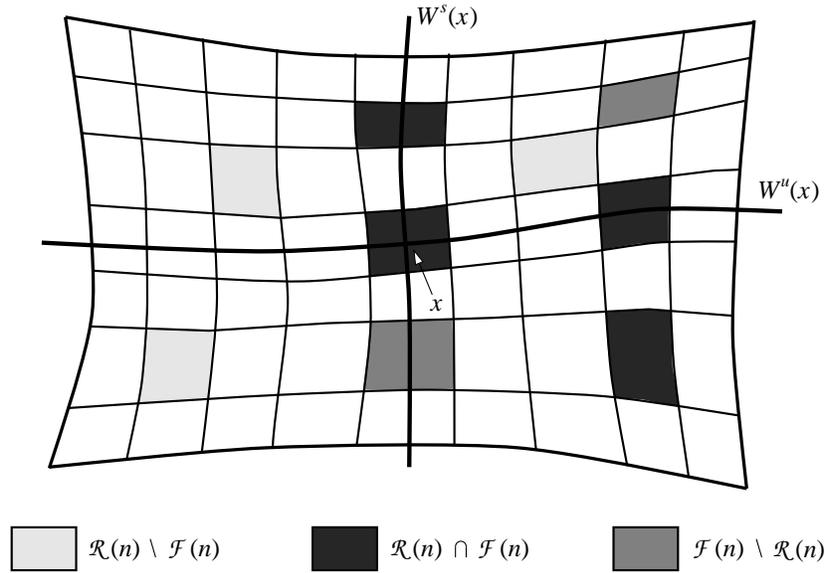

Figure 1. The procedure of "filling in" rectangles.

To effect this, for each set $A \subset \mathcal{P}(x)$, we define

$$N(n, A) = \text{card}\,\{R \in \mathcal{R}(n) : R \cap A \neq \emptyset\},$$
$$N^s(n, y, A) = \text{card}\,\{R \in \mathcal{R}(n) : R \cap \xi^s(y) \cap \Gamma \cap A \neq \emptyset\},$$
$$N^u(n, y, A) = \text{card}\,\{R \in \mathcal{R}(n) : R \cap \xi^u(y) \cap \Gamma \cap A \neq \emptyset\},$$
$$\widehat{N}^s(n, y, A) = \text{card}\,\{R \in \mathcal{F}(n) : R \cap \xi^s(y) \cap A \neq \emptyset\},$$
$$\widehat{N}^u(n, y, A) = \text{card}\,\{R \in \mathcal{F}(n) : R \cap \xi^u(y) \cap A \neq \emptyset\}.$$

Note that $N(n, \mathcal{P}(x))$ is the cardinality of the set $\mathcal{R}(n)$, and $N^s(n, y, \mathcal{P}(x))$ (respectively $N^u(n, y, \mathcal{P}(x))$) is the number of rectangles in $\mathcal{R}(n)$ that intersect $\Gamma$ and the stable (respectively unstable) local manifold at $y$. The product $\widehat{N}^s(n, y, \mathcal{P}(x)) \times \widehat{N}^u(n, y, \mathcal{P}(x))$ is the cardinality of the set $\mathcal{F}(n)$ for a "typical" point $y \in \mathcal{P}(x)$.

Let $Q_n(y)$ be the set defined by (15).

LEMMA 1. *For each $y \in \mathcal{P}(x) \cap \Gamma$ and integer $n \geq n_0$, we have:*

$$N^s(n, y, Q_n(y)) \leq \mu_y^s(B^s(y, 4e^{-n})) \cdot Ce^{anh+an\varepsilon};$$
$$N^u(n, y, Q_n(y)) \leq \mu_y^u(B^u(y, 4e^{-n})) \cdot Ce^{anh+an\varepsilon}.$$

*Proof.* It follows from (17) that

$$\mu_y^s(B^s(y, 4e^{-n})) \geq \mu_y^s(Q_n(y)) \geq N^s(n, y, Q_n(y))$$
$$\cdot \min\{\mu_y^s(R) : R \in \mathcal{R}(n) \text{ and } R \cap \xi^s(y) \cap Q_n(y) \cap \Gamma \neq \emptyset\}.$$

(Note that the condition $R \cap \xi^s(y) \cap Q_n(y) \neq \emptyset$ implies that $R \in \mathcal{R}(n)$.) Let $z \in R \cap \xi^s(y) \cap Q_n(y) \cap \Gamma$ for some $R \in \mathcal{R}(n)$. By Proposition 5 we obtain $\mu_y^s(R) = \mu_y^s(\mathcal{P}_{an}^0(z)) = \mu_z^s(\mathcal{P}_{an}^0(z))$. The first inequality in the lemma follows now from (6). The proof of the second inequality is similar. □

LEMMA 2. *For each $y \in \mathcal{P}(x) \cap \widehat{\Gamma}$ and integer $n \geq n_1$,*

$$\mu(B(y, e^{-n})) \leq N(n, Q_n(y)) \cdot 2Ce^{-2anh + 2an\varepsilon}.$$

*Proof.* It follows from (19) and (16) that

$$\frac{1}{2}\mu(B(y, e^{-n})) \leq \mu(B(y, e^{-n}) \cap \Gamma) \leq \mu(Q_n(y) \cap \Gamma)$$
$$\leq N(n, Q_n(y)) \cdot \max\{\mu(R) : R \in \mathcal{R}(n) \text{ and } R \cap Q_n(y) \neq \emptyset\}.$$

(Note that the condition $R \cap Q_n(y) \neq \emptyset$ implies that $R \in \mathcal{R}(n)$.) The desired inequality follows from (5). □

We now estimate the number of rectangles in the classes $\mathcal{R}(n)$ and $\mathcal{F}(n)$.

LEMMA 3. *For $\mu$-almost every $y \in \mathcal{P}(x) \cap \widehat{\Gamma}$ there is an integer $n_2(y) \geq n_1$ such that for each $n \geq n_2(y)$, we have:*

$$N(n+2, Q_{n+2}(y)) \leq \widehat{N}^s(n, y, Q_n(y)) \cdot \widehat{N}^u(n, y, Q_n(y)) \cdot 2C^2 e^{4a(h+\varepsilon)} e^{4an\varepsilon}.$$

*Proof.* By the Borel density lemma (with $A = \widehat{\Gamma}$), for $\mu$-almost every $y \in \widehat{\Gamma}$ there is an integer $n_2(y) \geq n_1$ such that for all $n \geq n_2(y)$,

$$2\mu(B(y, e^{-n}) \cap \widehat{\Gamma}) \geq \mu(B(y, e^{-n})).$$

Since $\widehat{\Gamma} \subset \Gamma$, it follows from (16) that for all $n \geq n_2(y)$,

(22) $\quad 2\mu(Q_n(y) \cap \widehat{\Gamma}) \geq 2\mu(B(y, e^{-n}) \cap \widehat{\Gamma}) \geq \mu(B(y, e^{-n}))$
$$\geq \mu(B(y, 4e^{-n-2})) \geq \mu(Q_{n+2}(y)).$$

For any $m \geq n_2(y)$, by (5) and property (e), we have

$$\mu(Q_m(y)) = \sum_{\mathcal{P}_{am}^{am}(z) \subset Q_m(y)} \mu(\mathcal{P}_{am}^{am}(z)) \geq N(m, Q_m(y)) \cdot C^{-1} e^{-2amh - 2am\varepsilon}.$$

Similarly, for every $n \geq n_2(y)$, we obtain

$$\mu(Q_n(y) \cap \widehat{\Gamma}) = \sum_{\mathcal{P}_{an}^{an}(z) \subset Q_n(y)} \mu(\mathcal{P}_{an}^{an}(z) \cap \widehat{\Gamma}) \leq N_n \cdot Ce^{-2anh + 2an\varepsilon},$$

where $N_n$ is the number of rectangles $\mathcal{P}_{an}^{an}(z) \in \mathcal{R}(n)$ that have nonempty intersection with $\widehat{\Gamma}$.

Set $m = n + 2$. The last two inequalities together with (22) imply that

$$(23) \qquad N(n+2, Q_{n+2}(y)) \leq N_n \cdot 2C^2 e^{4a(h+\varepsilon)+4an\varepsilon}.$$

On the other hand, since $y \in \widehat{\Gamma}$ the intersections $\mathcal{P}_0^{an}(y) \cap \xi^u(y) \cap \widehat{\Gamma}$ and $\mathcal{P}_{an}^0(y) \cap \xi^s(y) \cap \widehat{\Gamma}$ are nonempty.

Consider a rectangle $\mathcal{P}_{an}^{an}(v) \subset Q_n(y)$ that has nonempty intersection with $\widehat{\Gamma}$. Then the rectangles $\mathcal{P}_{an}^0(v) \cap \mathcal{P}_0^{an}(y)$ and $\mathcal{P}_{an}^0(y) \cap \mathcal{P}_0^{an}(v)$ are in $\mathcal{F}(n)$ and intersect respectively the stable and unstable local manifolds at $y$. Hence, to any rectangle $\mathcal{P}_{an}^{an}(v) \subset Q_n(y)$ with nonempty intersection with $\widehat{\Gamma}$, one can associate the pair of rectangles $(\mathcal{P}_{an}^0(v) \cap \mathcal{P}_0^{an}(y), \mathcal{P}_{an}^0(y) \cap \mathcal{P}_0^{an}(v))$ in

$$\{R \in \mathcal{F}(n) : R \cap \xi^s(y) \cap Q_n(y) \neq \emptyset\} \times \{R \in \mathcal{F}(n) : R \cap \xi^u(y) \cap Q_n(y) \neq \emptyset\}.$$

Clearly this correspondence is injective. Therefore,

$$\widehat{N}^s(n, y, Q_n(y)) \cdot \widehat{N}^u(n, y, Q_n(y)) \geq N_n.$$

The desired inequality follows from (23). $\square$

Our next goal is to compare the growth rate in $n$ of the number of rectangles in $\mathcal{F}(n)$ with the number of rectangles in $\mathcal{R}(n)$. We start with an auxiliary result.

LEMMA 4. *For each $x \in \widehat{\Gamma}$ and integer $n \geq n_1$, we have:*

$$\widehat{N}^s(n, x, \mathcal{P}(x)) \leq D^{-1} C^2 e^{anh+3an\varepsilon};$$
$$\widehat{N}^u(n, x, \mathcal{P}(x)) \leq D^{-1} C^2 e^{anh+3an\varepsilon}.$$

*Proof.* Since the partition $\mathcal{P}$ is countable we can find points $y_i$ such that the union of the rectangles $\mathcal{P}_0^{an}(y_i)$ is $\mathcal{P}(x)$, and these rectangles are mutually disjoint. Without loss of generality we can assume that $y_i \in \widehat{\Gamma}$ whenever $\mathcal{P}_0^{an}(y_i) \cap \widehat{\Gamma} \neq \emptyset$. We have

$$(24) \quad N(n, \mathcal{P}(x)) \geq \sum_i N^s(n, y_i, \mathcal{P}_0^{an}(y_i)) \geq \sum_{i : \mathcal{P}_0^{an}(y_i) \cap \widehat{\Gamma} \neq \emptyset} N^s(n, y_i, \mathcal{P}_0^{an}(y_i)).$$

We now estimate $N^s(n, y_i, \mathcal{P}_0^{an}(y_i))$ from below for $y_i \in \widehat{\Gamma}$. By Propositions 4 and 5, and (6),

$$(25) \quad N^s(n, y_i, \mathcal{P}_0^{an}(y_i)) \geq \frac{\mu_{y_i}^s(\mathcal{P}_0^{an}(y_i) \cap \Gamma)}{\max\{\mu_z^s(\mathcal{P}_{an}^{an}(z)) : z \in \xi^s(y_i) \cap \mathcal{P}(x) \cap \Gamma\}}$$

$$\geq \frac{D}{\max\{\mu_z^s(\mathcal{P}_{an}^{an}(z)) : z \in \xi^s(y_i) \cap \mathcal{P}(x) \cap \Gamma\}}$$

$$= \frac{D}{\max\{\mu_z^s(\mathcal{P}_{an}^0(z)) : z \in \xi^s(y_i) \cap \mathcal{P}(x) \cap \Gamma\}}$$

$$\geq DC^{-1} e^{anh - an\varepsilon}.$$

Similarly (5) implies that

$$(26) \quad N(n, \mathcal{P}(x)) \leq \frac{\mu(\mathcal{P}(x))}{\min\{\mu(\mathcal{P}_{an}^{an}(z)) : z \in \mathcal{P}(x) \cap \Gamma\}} \leq C e^{2anh + 2an\varepsilon}.$$

We now observe that

$$(27) \quad \widehat{N}^u(n, x, \mathcal{P}(x)) = \text{card}\{i : \mathcal{P}_0^{an}(y_i) \cap \widehat{\Gamma} \neq \emptyset\}.$$

Putting (24), (25), (26), and (27) together we conclude that

$$Ce^{2anh + 2an\varepsilon} \geq N(n, \mathcal{P}(x))$$

$$\geq \sum_{i : \mathcal{P}_0^{an}(y_i) \cap \widehat{\Gamma} \neq \emptyset} N^s(n, y_i, \mathcal{P}_0^{an}(y_i))$$

$$\geq \widehat{N}^u(n, x, \mathcal{P}(x)) \cdot DC^{-1} e^{anh - an\varepsilon}.$$

This yields $\widehat{N}^u(n, x, \mathcal{P}(x)) \leq D^{-1} C^2 e^{anh + 3an\varepsilon}$. The other inequality can be proved in a similar way. □

We emphasize that the procedure of "filling in" rectangles to obtain the class $\mathcal{F}(n)$ may substantially increase the number of rectangles in the neighborhood of some points. However, the next lemma shows that this procedure of "filling in" does not add too many rectangles at almost every point.

LEMMA 5. *For $\mu$-almost every $y \in \mathcal{P}(x) \cap \widehat{\Gamma}$, we have:*

$$\varlimsup_{n \to +\infty} \frac{\widehat{N}^s(n, y, Q_n(y))}{N^s(n, y, Q_n(y))} e^{-7an\varepsilon} < 1;$$

$$\varlimsup_{n \to +\infty} \frac{\widehat{N}^u(n, y, Q_n(y))}{N^u(n, y, Q_n(y))} e^{-7an\varepsilon} < 1.$$

*Proof.* By (17) and (20), for each $n \geq n_1$ and $y \in \widehat{\Gamma}$,

$$\mu_y^s(Q_n(y)) \geq \mu_y^s(B^s(y, e^{-n}) \cap \Gamma) \geq \frac{1}{2} \mu_y^s(B^s(y, e^{-n})).$$

Since $\mathcal{P}_{an}^{an}(z) \subset \mathcal{P}_{an}^0(z)$ for every $z$, by virtue of (6) and (10) we obtain

$$
(28) \quad N^s(n, y, Q_n(y)) \geq \frac{\mu_y^s(Q_n(y))}{\max\{\mu_z^s(\mathcal{P}_{an}^{an}(z)) : z \in \xi^s(y) \cap \mathcal{P}(x) \cap \Gamma\}}
$$

$$
\geq \frac{1}{2} \frac{\mu_y^s(B^s(y, e^{-n}))}{\max\{\mu_z^s(\mathcal{P}_{an}^0(z)) : z \in \xi^s(y) \cap \mathcal{P}(x) \cap \Gamma\}}
$$

$$
\geq \frac{1}{2C} \frac{e^{-d^s n - n\varepsilon}}{e^{-anh + an\varepsilon}}.
$$

Let us consider the set

$$
F = \left\{ y \in \widehat{\Gamma} : \varliminf_{n \to +\infty} \frac{\widehat{N}^s(n, y, Q_n(y))}{N^s(n, y, Q_n(y))} e^{-7an\varepsilon} \geq 1 \right\}.
$$

For each $y \in F$ there exists an increasing sequence $\{m_j\}_{j=1}^\infty = \{m_j(y)\}_{j=1}^\infty$ of positive integers such that

$$
(29) \quad \widehat{N}^s(m_j, y, Q_{m_j}(y)) \geq \frac{1}{2} N^s(m_j, y, Q_{m_j}(y)) e^{7am_j\varepsilon}
$$

$$
\geq \frac{1}{4C} e^{-d^s m_j + am_j h + 5am_j \varepsilon}
$$

for all $j$ (note that $a > 1$).

We wish to show that $\mu(F) = 0$. Assume on the contrary that $\mu(F) > 0$. Let $F' \subset F$ be the set of points $y \in F$ for which there exists the limit

$$
\lim_{r \to 0} \frac{\log \mu_y^s(B^s(y, r))}{\log r} = d^s.
$$

Clearly $\mu(F') = \mu(F) > 0$. Then we can find $y \in F$ such that

$$
\mu_y^s(F) = \mu_y^s(F') = \mu_y^s(F' \cap \mathcal{P}(y) \cap \xi^s(y)) > 0.
$$

It follows from Frostman's lemma that

$$
(30) \quad \dim_H(F' \cap \xi^s(y)) = d^s.
$$

Let us consider the countable collection of balls

$$
\mathfrak{B} = \{B(z, 4e^{-m_j(z)}) : z \in F' \cap \xi^s(y);\ j = 1, 2, \ldots\}.
$$

By the Besicovitch covering lemma (see, for example, [3]) one can find a subcover $\mathfrak{C} \subset \mathfrak{B}$ of $F' \cap \xi^s(y)$ of arbitrarily small diameter and finite multiplicity $\rho = \rho(\dim M)$. This means that for any $L > 0$ one can choose a sequence of points $\{z_i \in F' \cap \xi^s(y)\}_{i=1}^\infty$ and a sequence of integers $\{t_i\}_{i=1}^\infty$, where $t_i \in \{m_j(z_i)\}_{j=1}^\infty$ and $t_i > L$ for each $i$, such that the collection of balls

$$
\mathfrak{C} = \{B(z_i, 4e^{-t_i}) : i = 1, 2, \ldots\}
$$

comprises a cover of $F' \cap \xi^s(y)$ whose multiplicity does not exceed $\rho$. We write $Q(i) = Q_{t_i}(z_i)$. The Hausdorff sum corresponding to this cover is

$$\sum_{B \in \mathfrak{C}} (\operatorname{diam} B)^{d^s - \varepsilon} = (8^{d^s - \varepsilon}) \sum_{i=1}^{\infty} e^{-t_i(d^s - \varepsilon)}.$$

By (29), we obtain

$$\sum_{i=1}^{\infty} e^{-t_i(d^s - \varepsilon)} \leq \sum_{i=1}^{\infty} \widehat{N}^s(t_i, z_i, Q(i)) \cdot 4Ce^{-at_i h - 4at_i \varepsilon}$$

$$\leq 4C \sum_{q=1}^{\infty} e^{-aqh - 4aq\varepsilon} \sum_{i: t_i = q} \widehat{N}^s(q, z_i, Q(i)).$$

Since the multiplicity of the subcover $\mathfrak{C}$ is at most $\rho$, each set $Q(i)$ appears in the sum $\sum_{i: t_i = q} \widehat{N}^s(q, z_i, Q(i))$ at most $\rho$ times. Hence,

$$\sum_{i: t_i = q} \widehat{N}^s(q, z_i, Q(i)) \leq \rho \widehat{N}^s(q, y, \mathcal{P}(y)).$$

From Lemma 4 it follows that

$$\sum_{B \in \mathfrak{C}} (\operatorname{diam} B)^{d^s - \varepsilon} \leq 4(8^{d^s - \varepsilon}) C \sum_{q=1}^{\infty} e^{-aqh - 4aq\varepsilon} \rho \widehat{N}^s(q, y, \mathcal{P}(y))$$

$$\leq 4(8^{d^s - \varepsilon}) D^{-1} C^3 \rho \sum_{q=1}^{\infty} e^{-aqh - 4aq\varepsilon + aqh + 3aq\varepsilon}$$

$$= 4(8^{d^s - \varepsilon}) D^{-1} C^3 \rho \sum_{q=1}^{\infty} e^{-aq\varepsilon} < \infty.$$

Since $L$ can be chosen arbitrarily large (and so also the numbers $t_i$), it follows that $\dim_H(F' \cap \xi^s(y)) \leq d^s - \varepsilon < d^s$. This contradicts (30). Hence $\mu(F) = 0$ and this yields the first inequality in the lemma. The proof of the second inequality is similar. $\square$

By Lemma 5, for $\mu$-almost every $y \in \mathcal{P}(x) \cap \widehat{\Gamma}$ there exists an integer $n_3(y) \geq n_2(y)$ such that if $n \geq n_3(y)$, then

(31) $$\widehat{N}^s(n, y, Q_n(y)) < N^s(n, y, Q_n(y)) e^{7an\varepsilon},$$

(32) $$\widehat{N}^u(n, y, Q_n(y)) < N^u(n, y, Q_n(y)) e^{7an\varepsilon}.$$

Moreover, by Lusin's theorem, for every $\varepsilon > 0$ there exists a subset $\Gamma_\varepsilon \subset \widehat{\Gamma}$ such that

$$\mu(\Gamma_\varepsilon) > \mu(\widehat{\Gamma}) - \varepsilon \quad \text{and} \quad n_\varepsilon \stackrel{\text{def}}{=} \sup\{n_1, n_3(y) : y \in \Gamma_\varepsilon\} < \infty,$$

and the inequalities (31) and (32) hold for every $n \geq n_\varepsilon$ and $y \in \Gamma_\varepsilon$.

A version of the following statement was proved by Ledrappier and Young in [16]. For the sake of the reader's convenience we include its proof.

LEMMA 6. *For every $\varepsilon > 0$, if $y \in \Gamma_\varepsilon$ and $n \geq n_\varepsilon$, then*

$$\mu_y^s(B^s(y, e^{-n}))\mu_y^u(B^u(y, e^{-n})) \leq \mu(B(y, 4e^{-n})) \cdot 4C^3 e^{11an\varepsilon}.$$

*Proof.* Let $z \in \Gamma_\varepsilon \cap Q_n(y)$. By (32), if $n \geq n_\varepsilon$ then

$$N^u(n, y, Q_n(y)) \leq \widehat{N}^u(n, y, Q_n(y)) = \widehat{N}^u(n, z, Q_n(y)) < N^u(n, z, Q_n(y))e^{7an\varepsilon}$$

and

$$(33) \qquad N^u(n, y, Q_n(y)) \leq \inf\{N^u(n, z, Q_n(y)) : z \in \Gamma_\varepsilon \cap Q_n(y)\}e^{7an\varepsilon}.$$

Since $N(n, Q_n(y))$ is equal to the number of rectangles $R \subset Q_n(y)$ we have

$$\widehat{N}^s(n, y, Q_n(y)) \times \inf\{N^u(n, z, Q_n(y)) : z \in Q_n(y)\} \leq N(n, Q_n(y)).$$

By (33), if $y \in \Gamma_\varepsilon$ and $n \geq n_\varepsilon$ then

$$N^s(n, y, Q_n(y)) \times N^u(n, y, Q_n(y)) \leq N(n, Q_n(y))e^{7an\varepsilon}.$$

As in (28), if $y \in \Gamma_\varepsilon$ and $n \geq n_\varepsilon$ then

$$N^s(n, y, Q_n(y)) \geq \mu_y^s(B^s(y, e^{-n})) \cdot (2C)^{-1} e^{anh - an\varepsilon},$$

$$N^u(n, y, Q_n(y)) \geq \mu_y^u(B^u(y, e^{-n})) \cdot (2C)^{-1} e^{anh - an\varepsilon}.$$

Moreover, by (5) and (16),

$$N(n, Q_n(y)) \leq \frac{\mu(Q_n(y))}{\min\{\mu(\mathcal{P}_{an}^{an}(z)) : z \in Q_n(y) \cap \Gamma\}} \leq \mu(B(y, 4e^{-n})) \cdot C e^{2anh + 2an\varepsilon}.$$

Putting together these inequalities we obtain the desired statement. $\square$

## 6. Proof of the Main Theorem

Given $\varepsilon > 0$, let the set $\widehat{\Gamma}$ be as in the previous sections. By Lemmas 2 and 3, for $\mu$-almost every $y \in \mathcal{P}(x) \cap \widehat{\Gamma}$ and $n \geq n_2(y)$, we obtain

$$\mu(B(y, e^{-n-2})) \leq \widehat{N}^s(n, y, Q_n(y)) \cdot \widehat{N}^u(n, y, Q_n(y)) \cdot 4C^3 e^{4a(h+\varepsilon)} e^{-2anh + 6an\varepsilon}.$$

By (31), (32), and Lemma 1 we obtain

(34)
$$\mu(B(y, e^{-n-2})) \leq N^s(n, y, Q_n(y)) \cdot N^u(n, y, Q_n(y)) \cdot 4C^3 e^{4a(h+\varepsilon)} e^{-2anh+20an\varepsilon}$$
$$\leq \mu_y^s(B^s(y, 4e^{-n}))\mu_y^u(B^u(y, 4e^{-n})) \cdot 4C^5 e^{4a(h+\varepsilon)} e^{22an\varepsilon}.$$

Statement 1 of the theorem follows now immediately from (34) and Lemma 6.

By (34) and statement 1 in Proposition 2, we obtain

$$\lim_{n \to +\infty} \frac{\log \mu(B(y, e^{-n}))}{-n} \geq d^s + d^u - 22a\varepsilon$$

for $\mu$-almost every $y \in \widehat{\Gamma}$. Since $\mu(\widehat{\Gamma}) > 1 - \varepsilon$ and $\varepsilon > 0$ is arbitrarily small, we conclude that

$$\underline{d} = \varliminf_{r \to 0} \frac{\log \mu(B(y, r))}{\log r} = \varliminf_{n \to +\infty} \frac{\log \mu(B(y, e^{-n}))}{-n} \geq d^s + d^u$$

for $\mu$-almost every $y \in M$. Combining this property with statement 2 in Proposition 2, we obtain statement 2 of the theorem. □

## 7. The case of nonergodic measures

We show how to modify our arguments in the case when the measure $\mu$ is not ergodic. Notice that in this case Proposition 2 is still valid, i.e., for $\mu$-almost every $x \in M$ the stable and unstable dimensions $d^s(x)$ and $d^u(x)$ exist but may depend on the point. Moreover, given $\varepsilon > 0$, there is a set $\Gamma'$ of measure $\geq 1 - \varepsilon$ and, for every $x \in \Gamma'$, a number $h(x)$ such that (5)–(18) hold if $h$ is replaced by $h(x)$, $d^s$ by $d^s(x)$ and $d^u$ by $d^u(x)$. Let us fix $\kappa > 0$ and consider the sets

$$\Gamma(x) = \{y \in M : |h(x) - h(y)| < \kappa, \ |d^s(x) - d^s(y)| < \kappa, \ |d^u(x) - d^u(y)| < \kappa\}.$$

The collection of these sets covers $\Gamma'$. Moreover, there exists a countable subcollection of sets $\{\Gamma^i\}_{i \in \mathbb{N}}$ which still covers $\Gamma'$. Let $\mu^i$ be the conditional measure generated by $\mu$ on $\Gamma^i$. We can apply the arguments in the proof of the Main Theorem to the measures $\mu^i$ and show that for almost every $x \in \Gamma^i$ the following inequality holds

$$d^s(x) + d^u(x) - c\kappa < d^i(x),$$

where $d^i(x)$ is the pointwise dimension of the measure $\mu^i$ and $c$ does not depend on $x$ or $\kappa$. Since the cover $\{\Gamma^i\}_{i \in \mathbb{N}}$ is countable letting $\kappa$ go to zero yields that for $\mu$-almost every $x \in M$

$$d^s(x) + d^u(x) \leq \underline{d}(x).$$

It follows from [16] that $\overline{d}(x) \leq d^s(x) + d^u(x)$.

## Appendix. Lipschitz property of intermediate foliations

Let $M$ be a smooth compact connected Riemannian manifold with distance $d$ and $f\colon M \to M$ a $C^{1+\alpha}$ diffeomorphism. We denote the Lyapunov exponents at a point $x \in M$ by

$$\lambda_1(x) < \cdots < \lambda_{k(x)}(x).$$

Recall that a point $x \in M$ is called *Lyapunov regular* if there exists a decomposition $T_x M = \bigoplus_{i=1}^{k(x)} E^i(x)$, where $E^i(x)$ is the vector space

$$E^i(x) = \left\{ v \in T_x M \setminus \{0\} : \lim_{n \to \pm\infty} \frac{1}{n} \log \|d_x f^n v\| = \lambda_i(x) \right\} \cup \{0\}$$

such that

$$\lim_{n \to \pm\infty} \frac{1}{n} \log |\det d_x f^n| = \sum_{i=1}^{k(x)} \lambda_i(x) \dim E^i(x).$$

Note that $\lambda_i(f(x)) = \lambda_i(x)$ and $d_x f E^i(x) = E^i(f(x))$ for each $i$.

We denote by $\Lambda \subset M$ the set of all Lyapunov regular points. Below we describe some notions and results on nonuniformly hyperbolic dynamical systems. For references, see [18], [14], [9], and [10].

Set $s(x) = \max\{i : \lambda_i(x) < 0\}$ and $u(x) = \min\{i : \lambda_i(x) > 0\}$. Fix $x \in \Lambda$ and $r(x) > 0$. For every $i = 1, \ldots, s(x)$ one defines the $i^{\text{th}}$ *stable leaf* at $x$ by

$$W^i(x) = \left\{ y \in B(x, r(x)) : \limsup_{n \to +\infty} \frac{1}{n} \log d(f^n(x), f^n(y)) < \lambda_i(x) \right\},$$

and for every $i = u(x), \ldots, k(x)$ the $i^{\text{th}}$ *unstable leaf* at $x$ by

$$W^i(x) = \left\{ y \in B(x, r(x)) : \limsup_{n \to -\infty} \frac{1}{|n|} \log d(f^n(x), f^n(y)) < -\lambda_i(x) \right\}.$$

Clearly, $W^i(x) \subset W^{i+1}(x)$ if $i < s(x)$ and $W^i(x) \supset W^{i+1}(x)$ if $i \geq u(x)$. We write $W^s(x) = W^{s(x)}(x)$ and $W^u(x) = W^{u(x)}(x)$ and call them, respectively, the *stable* and *unstable manifolds* at $x$. One can prove that if $r(x)$ is sufficiently small than $W^i(x)$ is a $C^{1+\alpha}$ immersed submanifold.

For every $i \leq s(x)$, set $D^i(x) = \bigoplus_{j=1}^{i} E^j(x)$, and for every $i \geq u(x)$, set $D^i(x) = \bigoplus_{j=i}^{k(x)} E^j(x)$. If $i$ is such that $\lambda_i(x) \neq 0$, then $T_x W^i(x) = D^i(x)$.

Recall that a function $A\colon \Lambda \to \mathbb{R}$ is called *tempered* (with respect to $f$) if for every $x \in \Lambda$ we have

$$\lim_{n \to \pm\infty} \frac{1}{n} \log \|A(f^n(x))\| = 0.$$

Given $r \in (0, r(x))$ we denote by $B^i(x, r) \subset W^i(x)$ the ball centered at $x$ of radius $r$ with respect to the induced distance on $W^i(x)$.

Given $\tau > 0$ and $x \in \Lambda$, let $(U_x, \varphi_x)$ be a *Lyapunov chart* at $x$. This means the following:

1. $\varphi_x: U_x \to \mathbb{R}^n$ is a local diffeomorphism with the property that the spaces $\mathbb{E}^i = \varphi_x(\exp_x E^i(x))$ form an orthogonal decomposition of $\mathbb{R}^n$;

2. the subspaces $\mathbb{D}^i = \varphi_x(\exp_x D^i(x))$ are independent of $x$;

3. if $i = 1, \ldots, k(x)$ and $v \in E^i(x)$ then

(A1) $\quad e^{\lambda_i(x)-\tau}\|\varphi_x(\exp_x v)\| \leq \|\varphi_{f(x)}(\exp_{f(x)} d_x f v)\| \leq e^{\lambda_i(x)+\tau}\|\varphi_x(\exp_x v)\|;$

4. there is a constant $K$ and a tempered function $A: \Lambda \to \mathbb{R}$ such that if $y, z \in U_x$ then

$$K\|\varphi_x y - \varphi_x z\| \leq d(y,z) \leq A(x)\|\varphi_x y - \varphi_x z\|;$$

5. there exists $\rho(x) \in (0, r(x))$ such that $B^i(x, \rho(x)) \subset W^i(x) \cap U_x$ for every $x \in \Lambda$ and $i = 1, \ldots, k(x)$ with $\lambda_i(x) \neq 0$. Moreover, for $1 \leq i \leq s(x)$, the manifolds $\varphi_x(W^i(x))$ are graphs of smooth functions $g^i: \mathbb{D}^i \to \mathbb{D}^{i+1}$ and for $u(x) \leq i \leq k(x)$, of smooth functions $g^i: \mathbb{D}^i \to \mathbb{D}^{i-1}$; the first derivatives of $g^i$ are bounded by $1/3$.

It follows that if $i \leq s(x)$ then

$$f(W^i(x) \cap U_x) \subset W^i(f(x)) \cap U_{f(x)},$$

and if $i \geq u(x)$ then

$$f^{-1}(W^i(x) \cap U_x) \subset W^i(f^{-1}(x)) \cap U_{f^{-1}(x)}.$$

Each set $\Lambda$ can be decomposed into sets $\Lambda_\beta$ in which the numbers $k(x)$, $\dim E^i(x)$, and $\lambda_i(x)$ are constant for each $i$. For every ergodic measure $\mu$ invariant under $f$ there exists a unique $\beta$ for which the set $\Lambda_\beta$ has full $\mu$-measure. From now on we restrict our consideration to a subset $\Lambda_\beta \subset \Lambda$ and set $k(x) = k$, $s(x) = s$, $u(x) = u$, and $\lambda_i(x) = \lambda_i$ for each $i$ and $x \in \Lambda_\beta$.

Given $\ell > 0$, consider a set

$$\Lambda'_{\beta\ell} = \{x \in \Lambda_\beta : \rho(x) > 1/\ell,\ A(x) < \ell,$$
$$\angle(E^i(x), \bigoplus_{j \neq i} E^j(x)) > 1/\ell,\ i = 1\ldots,k\}.$$

Let $\Lambda_{\beta\ell}$ be the closure of $\Lambda'_{\beta\ell}$. For each $x \in \Lambda'_{\beta\ell}$ there exists an invariant decomposition $T_x M = \bigoplus_{i=1}^{k(x)} E^i(x)$, invariant $i^{\text{th}}$ stable and unstable leaves $W^i(x)$, and Lyapunov chart $(U_x, \phi_x)$ at $x$ with the above properties. In particular, the functions $\rho(x)$ and $A(x)$ can be extended to $\Lambda'_{\beta\ell}$ such that $\rho(x) > 1/\ell$, $A(x) < \ell$, and $\angle(E^i(x), \bigoplus_{j \neq i} E^j(x)) > 1/\ell$ for $i = 1\ldots,k$. The set $\Lambda_{\beta\ell}$ is obviously compact. We also have that $\Lambda_{\beta\ell} \subset \Lambda_{\beta(\ell+1)}$ and $\Lambda_\beta = \bigcup_{\ell>0} \Lambda_{\beta\ell}$ (mod 0).

Let us fix $c > 0$, $\ell > 0$, $x \in \Lambda_{\beta\ell}$, and $y' \in \Lambda_{\beta\ell} \cap B^{i+1}(x, c/\ell)$. For each $i < s$, consider two local smooth manifolds $\Sigma_x$ and $\Sigma_{y'}$ in $W^{i+1}(x)$, containing $x$ and $y'$, respectively and transversal to $W^i(z)$ for all $z \in \Lambda_{\beta\ell} \cap B^{i+1}(x, c/\ell)$. The holonomy map
$$\Pi_i = \Pi_i(\Sigma_x, \Sigma_{y'}) \colon \Sigma_x \cap \Lambda_{\beta\ell} \cap B^{i+1}(x, c/\ell) \to \Sigma_{y'}$$
is defined by
$$\Pi_i(x') = W^i(x') \cap \Sigma_{y'}$$
with $x' \in \Sigma_x$. This map is well-defined if $c$ is sufficiently small ($c$ may depend on $\ell$ but does not depend on $x$ and $y$).

Our goal is to prove the following theorem.

THEOREM. *Let $f$ be a $C^{1+\alpha}$ diffeomorphism. For each $\ell > 0$, $i < s$, $x \in \Lambda_{\beta\ell}$, and $y' \in \Lambda_{\beta\ell} \cap B^i(x, c/\ell)$ the holonomy map $\Pi_i(\Sigma_x, \Sigma_{y'})$ is Lipschitz continuous with the Lipschitz constant depending only on $\beta$ and $\ell$.*

*Proof.* We need a series of lemmas. First we observe that the Lipschitz property does not depend on the choice of transversal sections.

LEMMA A1. *For any local smooth manifolds $\Sigma'_x$ and $\Sigma'_{y'}$ in $W^{i+1}(x)$ containing $x$ and $y'$ respectively and transversal to $W^i(z)$ for $z \in B^{i+1}(x, c/\ell) \cap \Lambda_{\beta\ell}$, the holonomy map $\Pi_i(\Sigma'_x, \Sigma'_{y'})$ is Lipschitz continuous if and only if the map $\Pi_i(\Sigma_x, \Sigma_y)$ is Lipschitz continuous (possibly with a different Lipschitz constant depending only on the angles between the tangent spaces to the transversals).*

*Proof.* Since the foliation $W^i$ is continuous on the compact set
$$B^{i+1}(x, c/\ell) \cap \Lambda_{\beta\ell}$$
and its leaves are at least $C^1$, the angle between any transversal section to $W^i(z)$ (in $W^{i+1}(x)$) and any leaf of $W^i$ is uniformly bounded away from 0. Therefore, if $z \in \Sigma_x$ and $z' = \Sigma'_x \cap W^i(z)$ then the ratio $d(x,z)/d(x,z')$ is bounded from below and above (uniformly in $z$) by constants (depending only on those angles). This implies the desired result. □

For each $z \in M$ and $n \in \mathbb{N}$, we will write $z_n = f^n(z)$. Given $x \in \Lambda_\beta$, consider a map $\Phi_x \colon B^{i+1}(x, \rho(x)) \to \mathbb{D}^i \oplus \mathbb{E}^{i+1}$ by
$$\Phi_x y = (\pi_1(\varphi_x y), \pi_2(\varphi_x y)),$$
where $\pi_1 \colon \mathbb{D}^{i+1} \to \mathbb{D}^i$ and $\pi_2 \colon \mathbb{D}^{i+1} \to \mathbb{E}^{i+1}$ are the orthogonal projections. The following two lemmas are immediate consequences of the properties of the Lyapunov charts.

LEMMA A2. *The map $\Phi_x$ is a diffeomorphism with all derivatives tempered. Moreover, there is a constant $C_1$ (depending only on $\ell$) such that for every $x \in \Lambda_{\beta\ell}$, $n \in \mathbb{N}$, and $y, z \in B^{i+1}(x, \rho(x))$ we have*

1. $\|\Phi_{x_n}y_n - \Phi_{x_n}z_n\| \leq C_1 e^{(\lambda_{i+1}+\tau)n}\|\Phi_x y - \Phi_x z\|$;

2. if $\|\pi_1(\Phi_x y - \Phi_x z)\|/\|\pi_2(\Phi_x y - \Phi_x z)\| < 3\ell$ then

$$\|\pi_1(\Phi_{x_n}y_n - \Phi_{x_n}z_n)\| \leq C_1 e^{(\lambda_{i+1}+\tau)n}\|\pi_1(\Phi_x y - \Phi_x z)\|,$$
$$\|\pi_2(\Phi_{x_n}y_n - \Phi_{x_n}z_n)\| \leq C_1 e^{(\lambda_{i+1}+\tau)n}\|\pi_2(\Phi_x y - \Phi_x z)\|.$$

For each $x \in \Lambda_\beta$ we introduce the map $F_x \colon \Phi_x(B^{i+1}(x, \rho(x))) \to \mathbb{D}^i \oplus \mathbb{E}^{i+1}$ by

$$F_x(v) = \Phi_{f(x)} f(\Phi_x^{-1}(v)).$$

We write

$$F_x^n = \begin{cases} F_{x_n} \circ \cdots \circ F_x & n > 0 \\ \mathrm{Id} & n = 0 \\ F_{x_n}^{-1} \circ \cdots \circ F_x^{-1} & n < 0 \end{cases}.$$

LEMMA A3. *For every $x \in \Lambda_\beta$, the map $F_x$ is a $C^{1+\alpha}$ diffeomorphism onto its image, the Hölder constant $d(x)$ of the differential $dF_x(0)$ at $0$ is a tempered function of $x$, for every $\ell > 0$ there exists a constant $d_\ell > 0$ such that $d(x) \leq d_\ell$ for every $x \in \Lambda_{\beta\ell}$, and, finally, $dF_x(0)\mathbb{E}^i = \mathbb{E}^i$ for every $x \in \Lambda$ and every $i$.*

LEMMA A4. *There is a constant $C_2$ (depending only on $\ell$) such that if $n \in \mathbb{N}$, $x \in \Lambda_{\beta\ell}$, $y, z \in B^{i+1}(x, \rho(x))$, and*

$$\|\pi_1(\Phi_x y - \Phi_x z)\|/\|\pi_2(\Phi_x y - \Phi_x z)\| < 3\ell$$

*then*

(A2) $$C_2^{-1}\|\pi_2(\Phi_x y - \Phi_x z)\| \leq \|\pi_2 dF_{x_n}^{-n}(0)(\Phi_{x_n}y_n - \Phi_{x_n}z_n)\|$$
$$\leq C_2\|\pi_2(\Phi_x y - \Phi_x z)\|.$$

*Proof.* We first observe that

$$\|\Phi_{x_n}y_n - \Phi_{x_n}z_n\|$$
$$\leq \max\{\|\pi_1(\Phi_{x_n}y_n - \Phi_{x_n}z_n)\|, \|\pi_2(\Phi_{x_n}y_n - \Phi_{x_n}z_n)\|\}$$
$$\leq \max\{C_1 e^{(\lambda_i+\tau)n}\|\pi_1(\Phi_x y - \Phi_x z)\|, C_1 e^{(\lambda_{i+1}+\tau)n}\|\pi_2(\Phi_x y - \Phi_x z)\|\}$$
$$\leq \max\{C_1 e^{(\lambda_i+\tau)n}3\ell\|\pi_2(\Phi_x y - \Phi_x z)\|, C_1 e^{(\lambda_{i+1}+\tau)n}\|\pi_2(\Phi_x y - \Phi_x z)\|\}$$
$$\leq C_3 e^{(\lambda_{i+1}+\tau)n}\|\pi_2(\Phi_x y - \Phi_x z)\|$$
$$\leq C_4 e^{2\tau n}\|\pi_2(\Phi_{x_n}y_n - \Phi_{x_n}z_n)\|,$$

where $C_3 > 0$ and $C_4 > 0$ are constants. By Lemma A3 there is a constant $C_5 > 0$ (depending only on $\ell$) such that if $n \in \mathbb{N}$ then

$$\|\pi_2(\Phi_{x_{n-1}} y_{n-1} - \Phi_{x_{n-1}} z_{n-1} - dF_{x_n}^{-1}(\Phi_{x_n} y_n)(\Phi_{x_n} y_n - \Phi_{x_n} z_n))\|$$
$$\leq \|\Phi_{x_{n-1}} y_{n-1} - \Phi_{x_{n-1}} z_{n-1} - dF_{x_n}^{-1}(\Phi_{x_n} y_n)(\Phi_{x_n} y_n - \Phi_{x_n} z_n)\|$$
$$\leq C_5 e^{\tau n} \|\Phi_{x_n} y_n - \Phi_{x_n} z_n\|^{1+\alpha}.$$

If $\pi_2 \Phi_x y = \pi_2 \Phi_x z$, then the statement follows from the Lyapunov block form of $dF_{x_n}^{-n}(0)$. Assume now that $\pi_2 \Phi_x y \neq \pi_2 \Phi_x z$. By Lemma A2 and (A1) there are positive constants $C_6$ and $C_7$ (depending only on $\ell$) such that if $n \in \mathbb{N}$ then, letting $w_n = \Phi_{x_n} y_n - \Phi_{x_n} z_n$,

$$\frac{\|\pi_2 dF_{x_n}^{-n}(0)(\Phi_{x_n} y_n - \Phi_{x_n} z_n)\|}{\|\pi_2 dF_{x_{n-1}}^{-(n-1)}(0)(\Phi_{x_{n-1}} y_{n-1} - \Phi_{x_{n-1}} z_{n-1})\|}$$
$$\leq 1 + \frac{\|\pi_2 dF_{x_{n-1}}^{-(n-1)}(0)(w_{n-1} - dF_{x_n}^{-1}(0) w_n)\|}{\|\pi_2 dF_{x_{n-1}}^{-(n-1)}(0) w_{n-1}\|}$$
$$\leq 1 + \frac{\|dF_{x_{n-1}}^{-(n-1)}(0)[\pi_2(w_{n-1} - dF_{x_n}^{-1}(0) w_n)]\|}{\|dF_{x_{n-1}}^{-(n-1)}(0)[\pi_2 w_{n-1}]\|}$$
$$\leq 1 + C_6 e^{\tau n} \frac{\|\pi_2(w_{n-1} - dF_{x_n}^{-1}(0) w_n)\|}{\|\pi_2 w_{n-1}\|}$$
$$\leq 1 + C_6 e^{\tau n} \frac{\|\pi_2(w_{n-1} - dF_{x_n}^{-1}(\Phi_{x_n} y_n) w_n)\|}{\|\pi_2 w_{n-1}\|}$$
$$\quad + C_6 e^{\tau n} \frac{\|\pi_2((dF_{x_n}^{-1}(\Phi_{x_n} y_n) - dF_{x_n}^{-1}(0)) w_n)\|}{\|\pi_2 w_{n-1}\|}$$
$$\leq 1 + C_6 e^{\tau n} \frac{c_1 e^{\tau n} \|w_n\|^{1+\alpha}}{\|\pi_2 w_{n-1}\|}$$
$$\quad + C_6 e^{\tau n} \frac{\|dF_{x_n}^{-1}(\Phi_{x_n} y_n) - dF_{x_n}^{-1}(0)\| \cdot \|w_n\|}{\|\pi_2 w_{n-1}\|}$$
$$\leq 1 + C_5 C_6 C_4 e^{4\tau n} \|w_n\|^\alpha + C_6 e^{\tau n} C_7 e^{\tau n} \|\Phi_{x_n} y_n\|^\alpha \cdot C_4 e^{2\tau n}$$
$$\leq 1 + (C_5 + C_7) C_6 C_4 C_1^\alpha e^{[(\lambda_{i+1} + \tau)\alpha + 4\tau] n}.$$

Set $C_8 = (C_5 + C_7)C_6C_3C_1^\alpha$. Choosing $\tau$ sufficiently small we may assume that $4\tau + (\lambda_{i+1} + \tau)\alpha \leq \gamma < 0$. Therefore,

$$\frac{\|\pi_2 dF_{x_n}^{-n}(0)(\Phi_{x_n}y_n - \Phi_{x_n}z_n)\|}{\|\pi_2(\Phi_x y - \Phi_x z)\|} \leq \prod_{n=0}^{\infty}\left(1 + C_8 e^{[4\tau+(\lambda_{i+1}+\tau)\alpha]n}\right)$$

$$\leq \prod_{n=0}^{\infty}(1 + C_8 e^{\gamma n}) = C_2 < +\infty$$

and the upper estimate in (A2) follows. The lower estimate can be obtained in a similar way. □

We now proceed with the proof of the theorem. By Lemma A1 we may specify the choice of the transversals.

We say that a transversal section $\Sigma$ is *good* if $\varphi_x(\Sigma)$ is the graph of a function

$$h^i \colon \varphi_x(B^{i+1}(x, \rho(x))) \cap \mathbb{E}^{i+1} \to \mathbb{D}^{i+1}$$

whose first derivative is bounded by $1/3$. If $2\tau < |\lambda_i - \lambda_{i+1}|$ and $x \in \Lambda_{\beta\ell}$, then, by the properties of the Lyapunov charts, all the images $f^n(\Sigma)$ are good provided $\Sigma$ is good. Moreover, there is a constant $C_0 > 0$ such that

(A3) $\qquad C_0^{-1} d_{\varphi_x \Sigma}(\varphi_x y, \varphi_x z) \leq \|\pi_2(\varphi_x y - \varphi_x z)\| \leq C_0 d_{\varphi_x \Sigma}(\varphi_x y, \varphi_x z)$

for any good transversal section where $d_{\varphi_x \Sigma}$ is the induced distance on $\varphi_x \Sigma$.

Let $x \in \Lambda_{\beta\ell}$. We fix good transversals $\Sigma_x$ and $\Sigma_{y'}$ with $y' \in B^{i+1}(x, c/\ell) \cap \Lambda_{\beta\ell}$. Let $y = \Pi_i(x)$ and $y' = \Pi_i(x')$. Since the transversals are good we know that

$$\|\pi_1(\Phi_x y - \Phi_x y')\|/\|\pi_2(\Phi_x y - \Phi_x y')\| < 3\ell.$$

We also have that

$$\|\pi_1(\Phi_x x')\|/\|\pi_2(\Phi_x x')\| < 3\ell.$$

Set $C = 10\ell^2 C_0{}^2 C_2{}^2$. By Lemma A4 and (A3), for every $n \in \mathbb{N}$

$$\frac{d(y,y')}{d(x,x')} \leq 10\frac{d_{\Sigma_y}(y,y')}{d_{\Sigma_x}(x,x')}$$

$$\leq 10\ell^2 \frac{d_{\varphi_x\Sigma_y}(\varphi_x y, \varphi_x y')}{d_{\varphi_x\Sigma_x}(0, \varphi_x x')} \leq 10\ell^2 C_0{}^2 \frac{\|\pi_2(\Phi_x y - \Phi_x y')\|}{\|\pi_2(\Phi_x x')\|}$$

$$\leq C\frac{\|\pi_2 dF_{x_n}^{-n}(0)(\Phi_{x_n} y_n - \Phi_{x_n} y_n')\|}{\|\pi_2 dF_{x_n}^{-n}(0)\Phi_{x_n} x_n'\|}$$

$$\leq C\left(1 + \frac{\|\pi_2 dF_{x_n}^{-n}(0)(\Phi_{x_n} x_n' - \Phi_{x_n} y_n' + \Phi_{x_n} y_n)\|}{\|\pi_2 dF_{x_n}^{-n}(0)\Phi_{x_n} x_n'\|}\right)$$

$$\leq C\left(1 + \frac{\|\pi_2 dF_{x_n}^{-n}(0)(\Phi_{x_n} x_n' - \Phi_{x_n} y_n')\| + \|\pi_2 dF_{x_n}^{-n}(0)\Phi_{x_n} y_n\|}{\|\pi_2 dF_{x_n}^{-n}(0)\Phi_{x_n} x_n'\|}\right)$$

$$\leq C\left(1 + \frac{e^{\tau n}d(x_n', y_n') + e^{\tau n}d(x_n, y_n)}{e^{(\lambda_{i+1} - \tau)n}}\right)$$

$$\leq C(1 + e^{(3\tau + \lambda_i - \lambda_{i+1})n}).$$

Choosing $\tau$ such that $3\tau + \lambda_i - \lambda_{i+1} < 0$, we conclude that $\Pi_i$ is Lipschitz. □

Instituto Superior Técnico, Lisboa, Portugal
*E-mail address*: barreira@math.ist.utl.pt
*URL*: http://www.math.ist.utl.pt/˜barreira/

The Pennsylvania State University, University Park, PA
*E-mail address*:pesin@math.psu.edu
*URL*: http://www.math.psu.edu/pesin/

Weierstrass Institute of Applied Analysis and Stochastics, Berlin, Germany
*Current address*: Freie Universität Berlin, Berlin, Germany
*E-mail*: shmeling@math.fu-berlin.de